\documentclass{amsart}
\newtheorem{thm}{Theorem}
\newtheorem{lem}[thm]{Lemma}

\theoremstyle{definition}

\begin{document}
\title[On the Chv\a'atal-Janson Conjecture]{On the Chv\a'atal-Janson Conjecture}

\author{Lucio Barabesi}
\address{Lucio Barabesi, Dipartimento di Economia Politica e Statistica, Universit\`a di Siena, via San Francesco 7, 53100 Siena, Italy}
\email{lucio.barabesi@unisi.it}

\author{Luca Pratelli}
\address{Luca Pratelli, Accademia Navale, viale Italia 72, 57100 Livorno, Italy}
\email{pratel@mail.dm.unipi.it}

\author{Pietro Rigo}
\address{Pietro Rigo (corresponding author), Dipartimento di Scienze Statistiche ''P. Fortunati'', Universita' di Bologna, via delle Belle Arti 41, 40126 Bologna, Italy}
\email{pietro.rigo@unibo.it}



\begin{abstract}
In a recent paper, Svante Janson has considered a conjecture
suggested by Va\v{s}ek Chv\a'atal dealing with the probability that a
binomial random variable with parameters $n$ and $m/n$ - where $m$ is an
integer - exceeds its expectation $m$. Albeit Janson has provided a proof of
this conjecture for large $n$, we show that the result actually holds
for each $n\geq 2$.
\end{abstract}

\maketitle

\section{Introduction}\label{intro}

 By assuming that $B(n,m/n)$ denotes a binomial random variable
with parameters $n$ and $m/n$, Janson (2021) introduces the following
conjecture suggested by Va\v{s}ek Chv\a'atal in a personal communication.
\\

\noindent{\bf Conjecture 1} (Chv\a'atal). {\it For any fixed $n\geq 2$, as $m$ ranges over
$\{0,\ldots,n\}$, the probability $q_m:=P(B(n,m/n)\leq m)$ is the smallest
when $m=[\![2n/3]\!]$ where $[\![\cdot]\!]$ represents the nearest integer
function.}
\\

It is worth noting that the conjecture may have interesting applications,
since the probability that a binomial random variable exceeds its expected
value has generally an important role in the machine learning literature (see
e.g. Doerr, 2018, Greenberg and Mohri, 2014, Vapnik, 1998). Such a
probability has even a connection with an equation given by Ramanujan, as
emphasized by Jogdeo and Samuels (1968). For further results on the topic,
see Pelekis and Ramon (2016), Slud (1977).

Janson (2021) has proven that, for large $n$, Conjecture 1 actually holds and
the probabilities $q_m$ have a unique minimum. More precisely, Janson (2021)
provides the following theorem.
\\

\noindent{\bf Theorem 1}. {\it There exists a $n_0$ such that for each $n_0\geq n$:
$i$) $q_m$ is minimum for $m=[\![2n/3]\!]$ and $ii)$ $q_m\geq q_{m+1}$ if and only if $m+\frac{1}{2}<2n/3$.}
\\

\noindent However, Janson (2021) remarks that, even if it is possible in principle
computing an explicit value for $n_0$ in the proof of Theorem 1 and
numerically checking the statement for $n<n_0$, such an issue is not
practically feasible. Actually, Janson (2021, Remark 1.5) wishes for a
general proof of Theorem 1.

In the present contribution, we give a plain proof of Theorem 1 for each
$n\geq 2$. The proof is achieved by means of different methods with respect
to those adopted by Janson (2021), which are actually based on the version
for integer-valued random variables of the asymptotic Edgeworth expansion for
probabilities in the central limit theorem - as proposed by Esseen (1945).
Indeed, our proof shares similarities with the approach introduced by
Rigollet and Tong (2011, Appendix B) for assessing that $q_m\geq q_{m+1}$ for
$0\leq m<n/2$ and $n\geq 2$.

\section{Notations and Preliminaries}\label{m9v4} Let $(U_i)_{1\leq i\leq n}$ be $n$ independent copies of a Uniform random variable on $[0,1]$. If $(U_{(i)})_{1\leq i\leq n}$ represent the order statistics corresponding to $(U_i)_{1\leq i\leq n}$, it obviously holds $$q_m=P(\sum_{i=1}^n I_{\{U_i\leq m/n\}}\leq m)=P(U_{(m+1)}>m/n).\eqno(1)$$
On the basis of (1), for each $m\leq n-1$ it follows that $$q_m=(m+1){n\choose {m+1}}\int_{m/n}^1x^m(1-x)^{n-m-1}\, dx,\eqno(2)$$
since the probability density function $f_{m+1}$ of $U_{(m+1)}$ is given by $$f_{m+1}(x)=(m+1){n\choose {m+1}}x^m(1-x)^{n-m-1}I_{[0,1]}(x)$$ (see e.g. Feller, 1971, Section I.7).

\medskip\medskip\medskip\noindent
\begin{lem}
Let us assume that $n\geq 2$ and $m\leq n-2$. Then $$q_m\geq q_{m+1}\Longleftrightarrow\int_{m/n}^{(m+1)/n}x^{m+1}(1-x)^{n-m-2}\, dx\geq b_m,$$ where $b_m={{(m/n)^{m+1}(1-m/n)^{n-m-1}}\over{n-m-1}}.$ In addition, $q_m\geq q_{m+1}$ is equivalent to
$$\int_{0}^{1}(1+t/m)^{m+1}(1-t/(n-m))^{n-m-2}\, dt
\geq {{n-m}\over{n-m-1}},\eqno(3)  $$ or to $$ \int_{0}^{1}(1-{{v}\over{m+1}})^{m}(1+{{v}\over{n-m-1}})^{n-m-1}\, dv
\geq 1.\eqno(4)$$
\end{lem}

\begin{proof} On the basis of $(2)$ and by using the definition of the binomial coefficient, it follows that $q_m\geq q_{m+1}$ is equivalent to
$${{m+1}\over{n-m-1}}\int_{m/n}^1x^m(1-x)^{n-m-1}\, dx\geq\int_{(m+1)/n}^1x^{m+1}(1-x)^{n-m-2}\, dx.$$
Integrating by part the left-hand side of the previous inequality, the expression reduces to
$$-b_m+\int_{m/n}^1x^{m+1}(1-x)^{n-m-2}\, dx\geq\int_{(m+1)/n}^1x^{m+1}(1-x)^{n-m-2}\, dx$$
and the main result follows. As to $(3)$, from the main result and by means of the substitution $x=m/n+t/n$, it reads
$$\int_{0}^{1}(m/n+t/n)^{m+1}(1-m/n-t/n)^{n-m-2}\, dt
\geq n b_m,$$ which provides $(3)$ by suitably dividing both sides by $(m/n)^{m+1}(1-m/n)^{n-m-2}$.
As to $(4)$, by multiplying both sides of $(3)$ by $(n-m-1)/(n-m)$ and integrating by parts the corresponding left-hand side, it reads
$$\int_{0}^{1}(1+{{t}\over{m}})^{m}(1-{{t}\over{n-m}})^{n-m-1}\, dt\geq (1-{{1}\over{n-m}})^{n-m-1}(1+{{1}\over{m}})^{m}.$$ By dividing both sides of the previous inequality by the quantity in the right-hand side, the  expression reduces to
$$\int_{0}^{1}({{t+m}\over{1+m}})^{m}({{n-m-t}\over{n-m-1}})^{n-m-1}\, dt\geq 1, $$ which readily provides $(4)$ by considering the transformation $t=1-v$.
\end{proof}

\medskip\medskip\smallskip
\begin{lem} For a given $n\geq3$, let $m$ be an integer in $[1,n-2]$. For each $v\in]0,1]$, the function $g_v$ defined on $[1,n-2]$ and such that $$g_v(m)=(1-{{v}\over{m+1}})^{m}(1+{{v}\over{n-m-1}})^{n-m-1}$$ is decreasing. Moreover, the function $$h:\ m\mapsto \int_{0}^{1}(1-{{v}\over{m+1}})^{m}(1+{{v}\over{n-m-1}})^{n-m-1}\, dv$$
is decreasing on $[1,n-2]$.
\end{lem}

\begin{proof} For a given $v\in]0,1]$ and by denoting with $x$ a real number in $[1,n-2]$, it suffices to prove that $g^\prime_v(x)<0$. Since
$$g^\prime_v(x)=g_v(x)\, \big[\log(1-{{v}\over{x+1}})+ {{{{{vx}\over{(x+1)^2}}}}\over{ 1- {{v}\over{x+1}} }} -\log(1+{{v}\over{n-x-1}})+
{{{{{v}\over{n-x-1}}}}\over{ 1+ {{v}\over{n-x-1}} }}\big],$$it holds $$g^\prime_v(x)<0\Longleftrightarrow {{{{{vx}\over{(x+1)^2}}}}\over{ 1- {{v}\over{x+1}} }}+{{{{{v}\over{n-x-1}}}}\over{ 1+ {{v}\over{n-x-1}} }}<\log\big({{1+{{v}\over{n-x-1}}}\over{ 1-{{v}\over{x+1}}  }}\big).$$
In addition, since
\begin{equation*}
\begin{aligned}
\log\big({{1+{{v}\over{n-x-1}}}\over{ 1-{{v}\over{x+1}}  }}\big)&=\log\big[\big(1+{{v}\over{n-x-1}}\big)\big( 1+{{{{{v}\over{x+1}}}}\over{ 1- {{v}\over{x+1}} }}\big)\big]\cr &= \log\big(1+{{v}\over{n-x-1}}\big)+\log \big( 1+{{{{{v}\over{x+1}}}}\over{ 1- {{v}\over{x+1}} }}\big),
\end{aligned}
\end{equation*}

\noindent in order to prove that $g^\prime_v(x)<0$ it suffices to show that $${{{{{v}\over{n-x-1}}}}\over{ 1+ {{v}\over{n-x-1}} }}<\log\big(1+{{v}\over{n-x-1}}\big)\eqno(5)$$ and
$$\log \big( 1+{{{{{v}\over{x+1}}}}\over{ 1- {{v}\over{x+1}} }}\big)>{{{{{vx}\over{(x+1)^2}}}}\over{ 1- {{v}\over{x+1}} }}.\eqno(6)$$ By assuming that $c={{v}\over{n-x-1}}$, inequality $(5)$ follows from $\log(1+c)>c/(1+c)$, which holds for each $c>0$, while inequality $(6)$ is equivalent to $$\log \big( 1+{{v}\over{x+1-v}} \big)-{{v}\over{x+1-v}}+{{v}\over{(x+1)(x+ 1- v)}}>0$$ which, by assuming that $c=v/(x+1-v)$, reduces to
$$\log(1+c)-c+{{c^2}\over{v(c+1)}}>0.\eqno(7)$$Inequality $(7)$ holds since $$
\log(1+c)-c+{{c^2}\over{v(c+1)}}\geq \log(1+c)-c+{{c^2}\over{c+1}}=\log(1+c)-{{c}\over{c+1}}$$and Lemma is proven.
\end{proof}

\medskip\section{Proof of the Chv\'atal-Janson conjecture}

\medskip\noindent In this section we provide a proof of Theorem 1 for $n\geq 2$. On the basis of Lemma 1 and Lemma 2, for a given $n=3s+r$, where $s\geq 1$ and $r\in\{0,1,2\}$, it suffices to prove that for $r=0$ it holds
$$ \int_{0}^{1}(1-{{v}\over{2s}})^{2s-1}(1+{{v}\over{s}})^{s}\, dv\geq 1>\int_{0}^{1}(1-{{v}\over{2s+1}})^{2s}(1+{{v}\over{s-1}})^{s-1}\, dv,\eqno(8) $$ while for $r=1,2$ it holds $$
\int_{0}^{1}(1-{{v}\over{2s+1}})^{2s}(1+{{v}\over{s+r-1}})^{s+r-1}\, dv\geq 1\eqno(9)$$ and $$1>\int_{0}^{1}(1-{{v}\over{2s+2}})^{2s+1}(1+{{v}\over{s+r-2}})^{s+r-2}\, dv.\eqno(10)  $$
By considering the inequalities $(8)$ (i.e. when $r=0$) and by applying the Bernoulli inequality$(1+c)^s\geq 1+sc$ which holds for each $c>-1$, it follows that the first inequality in $(8)$ is true for each $s\geq 1$ since

\begin{equation*}
\begin{aligned}
\int_{0}^{1}(1-{{v}\over{2s}})^{2s-1}(1+{{v}\over{s}})^{s}\, dv &=\int_{0}^{1}[(1-{{v}\over{2s}})^{2}(1+{{v}\over{s}})]^{s}(1-{{v}\over{2s}})^{-1}\, dv\cr &=\int_{0}^{1}[1-{{3v^2}\over{4s^2}}+{{v^3}\over{4s^3}}]^{s}(1-{{v}\over{2s}})^{-1}\, dv\cr &\geq \int_0^1(1-{{3v^2}\over{4s}}+{{v^3}\over{4s^2}})(1-{{v}\over{2s}})^{-1}\, dv\cr &\geq \int_0^1(1-{{3v^2}\over{4s}}+{{v^3}\over{4s^2}})(1+{{v}\over{2s}}+{{v^2}\over{4s^2}})\, dv\cr &
=1+{{5}\over{96s^2}}-{{1}\over{80s^3}}+{{1}\over{96s^4}}>1 .
\end{aligned}
\end{equation*}

\noindent The first inequality in the previous expression follows from $(1-c)^{-1}>1+c+c^2$, which holds for each $c\in]0,1[$. In turn for $r=0$, as to the second inequality in (8) and by assuming that
$I_s=\int_{0}^{1}(1-{{v}\over{2s+1}})^{2s}(1+{{v}\over{s-1}})^{s-1}\, dv$, it reads

\begin{equation*}
\begin{aligned}
I_s &=\int_{0}^{1}[(1-{{v}\over{2s+1}})^{2s}(1+{{v}\over{s-1}})]^{s}(1+{{v}\over{s-1}})^{-1}\, dv\cr &
=\int_{0}^{1}\kern-0.9mm{{2s+1}\over{2s+1-v}}\exp\big((2s\kern-0.4mm+\kern-0.4mm1)\log(1-{{v}\over{2s+1}})\kern-0.4mm+\kern-0.4mm
(s-1)\log(1\kern-0.4mm+\kern-0.4mm{{v}\over{s\kern-0.4mm-\kern-0.4mm1}})\big)
 dv
\cr &=\int_{0}^{1}{{2s+1}\over{2s+1-v}}\exp\big(
\sum_{k=2}^\infty {{v^k}\over{k}}({{(-1)^{k-1}}\over{(s-1)^{k-1}}}-{{1}\over{(2s+1)^{k-1}}})\big)
\, dv\cr &
<\int_{0}^{1}{{2s+1}\over{2s+1-v}}\exp\big(
\sum_{k=2}^3 {{v^k}\over{k}}({{(-1)^{k-1}}\over{(s-1)^{k-1}}}-{{1}\over{(2s+1)^{k-1}}})\big)
\, dv,
\end{aligned}
\end{equation*}

\noindent since $\sum_{k=4}^\infty {{v^k}\over{k}}({{(-1)^{k-1}}\over{(s-1)^{k-1}}}-
{{1}\over{(2s+1)^{k-1}}})<0$. By adopting the notation $$\gamma(s,v)=\sum_{k=2}^3 {{v^k}\over{k}}({{(-1)^{k-1}}\over{(s-1)^{k-1}}}-{{1}\over{(2s+1)^{k-1}}}),$$ it should be remarked that $$\gamma(s,v)=-{{3v^2s}\over{2(s-1)(2s+1)}}+{{v^3s(s+2)}\over{(s-1)^2(2s+1)^2}}<0$$ for each $s\geq 2$ and $v\in ]0,1]$. Since $\exp(c)< 1+c+{{c^2}\over{2}}$ for $c<0$, it follows

\begin{equation*}
\begin{aligned}
I_s &<\int_{0}^{1}{{2s+1}\over{2s+1-v}}\exp\big(\gamma(s,v)\big)
\, dv\cr &<\int_{0}^{1}{{2s+1}\over{2s+1-v}}\big(1+\gamma(s,v)+{{\gamma(s,v)^2}\over{2}}\big)
\, dv.
\end{aligned}
\end{equation*}

\noindent Moreover, since
\begin{equation*}
\begin{aligned}
{{2s+1}\over{2s+1-v}} &={{1}\over{1-{{v}\over{2s+1}}}}= 1+{{v}\over{2s+1}}+{{v^2}\over{(2s+1)^2}}
{{1}\over{1-{{v}\over{2s+1}}}}\cr &\leq
1+{{v}\over{2s+1}}+{{5v^2}\over{4(2s+1)^2}},
\end{aligned}
\end{equation*}

\noindent it also follows
\begin{equation*}
\begin{aligned}
I_s &<\int_{0}^{1}\big(1+{{v}\over{2s+1}}+
{{5v^2}\over{4(2s+1)^2}}\big)\big(1+\gamma(s,v)+{{\gamma(s,v)^2}\over{2}}\big)
\, dv.
\end{aligned}
\end{equation*}
By computing the integral and by means of tedious algebraic manipulations, it holds
\begin{equation*}
\begin{aligned}
I_s &<1+
{{1+3s(-11s^3+(s+4)^2)}\over{(s-1)^4(2s+1)^6}}
+
{{s^6(29+7s-15s^2/2)}\over{(s-1)^4(2s+1)^6}}.
\end{aligned}
\end{equation*}

\noindent Since the numerators of the two fractions in the previous expressions are negative for $s\geq 3$, it holds that $I_s<1$ for $s\geq 3$. Finally, by direct computation
it also follows that $I_2<1$ and hence Theorem 1 holds true for $r=0$. As to $(9)$, i.e. when $r=1,2$, let us assume that $$J_r=\int_{0}^{1}(1-{{v}\over{2s+1}})^{2s}(1+{{v}\over{s+r-1}})^{s+r-1}\, dv.$$ By remarking that for each $s$ it holds $(1+{{v}\over{s}})^{s}\leq (1+{{v}\over{s+1}})^{s+1}$, it follows $J_1\leq J_2$.
Hence, in order to prove $(9)$ it suffices to show $J_1\geq 1$. It holds
\begin{equation*}
\begin{aligned}
J_1 &=\int_{0}^{1}(1-{{v}\over{2s+1}})^{2s}(1+{{v}\over{s}})^{s}\, dv
 \cr &= \int_{0}^{1}[(1-{{v}\over{2s+1}})^{2}(1+{{v}\over{s}})]^s \, dv \cr &=\int_{0}^{1}[(1-{{2v}\over{2s+1}}+{{v^2}\over{(2s+1)^2}})(1+{{v}\over{s}})]^s\, dv\cr &=\int_{0}^{1}[(1+{{v(1-2v)}\over{s(2s+1)}}+{{v^2}\over{(2s+1)^2}}(1+{{v}\over{s}})]^s\, dv.
\end{aligned}
\end{equation*}

\noindent By applying Bernoulli inequality, it reads
\begin{equation*}
\begin{aligned}
J_1&\geq\int_{0}^{1}(1+{{v(1-2v)}\over{2s+1}}+{{sv^2}\over{(2s+1)^2}}(1+{{v}\over{s}}) \, dv\cr &= 1-{{1}\over{6(2s+1)}}+{{s}\over{3(2s+1)^2}}+
{{1}\over{4(2s+1)^2}}=1+{{1}\over{12(2s+1)^2}},
\end{aligned}
\end{equation*}

\noindent which obviously implies $(9)$.
Finally, we prove $(10)$.
By adopting the notation $$H_r=\int_{0}^{1}(1-{{v}\over{2s+2}})^{2s+1}(1+{{v}\over{s+r-2}})^{s+r-2}\, dv, $$
since for each $v$ it holds $(1+{{v}\over{s-1}})^{s-1}\leq (1+{{v}\over{s}})^{s},$ it also follows that $H_1\leq H_2$ and hence it suffices to show the case $H_2<1$. To this aim, similarly to the the proof of the inequality $I_s<1$, it reads
\begin{equation*}
\begin{aligned}
H_2 &=\int_{0}^{1}(1-{{v}\over{2s+2}})^{2s+1}(1+{{v}\over{s}})^{s}\, dv\cr &
=\int_{0}^{1}\kern-0.9mm{{2s+2}\over{2s+2-v}}\exp\big((2s\kern-0.4mm+\kern-0.4mm 2)\log(1-{{v}\over{2s+2}})\kern-0.4mm+\kern-0.4mm
s\log(1\kern-0.4mm+\kern-0.4mm{{v}\over{s}})\big)
 dv
\cr &=\int_{0}^{1}{{2s+2}\over{2s+2-v}}\exp\big(
\sum_{k=2}^\infty {{v^k}\over{k}}({{(-1)^{k-1}}\over{s^{k-1}}}-{{1}\over{(2s+2)^{k-1}}})\big)
\, dv\cr &
<\int_{0}^{1}{{2s+2}\over{2s+2-v}}\exp\big(
\sum_{k=2}^3 {{v^k}\over{k}}({{(-1)^{k-1}}\over{s^{k-1}}}-{{1}\over{(2s+2)^{k-1}}})\big)
\, dv,
\end{aligned}
\end{equation*}
since $\sum_{k=4}^\infty {{v^k}\over{k}}({{(-1)^{k-1}}\over{s^{k-1}}}-
{{1}\over{(2s+2)^{k-1}}})<0$.
By assuming that $$\lambda(s,v)=\sum_{k=2}^3 {{v^k}\over{k}}({{(-1)^{k-1}}\over{s^{k-1}}}-{{1}\over{(2s+2)^{k-1}}}),$$ it should be noticed that $\lambda(s,v)<0$ for each $s\geq 2$ and $v\in ]0,1]$. By considering the inequality $\exp(c)< 1+c+{{c^2}\over{2}}$ for $c<0$, it follows
\begin{equation*}
\begin{aligned}
H_2 &<\int_{0}^{1}{{2s+2}\over{2s+2-v}}\exp\big(\lambda(s,v)\big)
\, dv\cr &<\int_{0}^{1}{{2s+2}\over{2s+2-v}}\big(1+\lambda(s,v)+{{\lambda(s,v)^2}\over{2}}\big)
\, dv\cr &<
\int_{0}^{1}\big(1+{{v}\over{2s+2}}+
{{6v^2}\over{5(2s+2)^2}}\big)\big(1+
\lambda(s,v)+{{\lambda(s,v)^2}\over{2}}\big)\, dv,
\end{aligned}
\end{equation*}

 \noindent since it holds ${{2s+2}\over{2s+2-v}}< 1+{{v}\over{2s+2}}+{{6v^2}\over{5(2s+2)^2}} $ for each $s\geq 2$ and $v\in]0,1[$. By evaluating the previous integral and by suitable algebraic manipulations, the following inequality holds for $s\geq 2$
$$H_2<1+{{-2s^5-9(s^4+s^3)+8s^2+22s+18}
\over{64s(s+1)^6}}.$$ The right-hand side of the previous inequality is obviously less than $1$. Moreover, a direct computation provides $H_2 <1$ for $s=1$. Therefore, Theorem 1 is proven.


\begin{thebibliography}{00}

\bibitem{DB}
Doerr Benjamin,  (2018)
An elementary analysis of the probability that a binomial random variable exceeds its expectation,

{\em Statist. Probab. Lett.}, 139,  67-74.

\bibitem{ES}
Esseen C.G., (1945) Fourier analysis of distribution functions. A
mathematical study of the Laplace-Gaussian law, {\em\ Acta Mathematica}, 77, 1-125.
\bibitem{FE}
Feller (1971)


\bibitem{GS}
Greenberg, Spencer, Mohri, Mehryar, (2014) Tight lower bound on the probability of a binomial exceeding its
expectation, {\em Statist. Probab. Lett.}, 86, 91-98,

\bibitem{SJ}
 Svante Janson, (2021) On the probability that a binomial variable is at most its
expectation, {\em Statist. Probab. Lett.},  171, 109020.

\bibitem{JSK}
Jogdeo, K. and Samuels, S.M., (1968) Monotone convergence of binomial
probabilities and a generalization of Ramanujan's equation, {\em Annals of
Mathematical Statistics}, 39, 1191-1195.

\bibitem{PR}
Pelekis C. and Ramon J., (2016) A lower bound on the probability that a
binomial random variable is exceeding its mean, {\em Statistics and Probability
Letters}, 119, 305-309.

\bibitem{TMG}
 Rigollet Philippe, Tong Xin, (2011) Neyman–pearson classification, convexity and stochastic constraints, {\em J. Mach. Learn. Res.}, 12, 2831-2855.

\bibitem{SL}
Slud, E.V., (1977) Distribution inequalities for the binomial law, {\em Annals of
Probability}, 5, 404-412.

\bibitem{VV}
Vapnik V.N., (1998) Statistical Learning Theory, Wiley, New York.


\end{thebibliography}
\end{document}